\documentstyle[12pt]{article}
\newcommand{\const}{\mathop{\rm const}\limits}

\newcommand{\card}{\mathop{\rm card}\limits}

\newcommand{\supp}{\mathop{\rm supp}\limits}

\newcommand{\Var}{\mathop{\rm Var}\limits}

\newcommand{\rank}{\mathop{\rm rank}\limits}

\begin{document}

\begin{center}

{\bf SHARP MOMENT AND EXPONENTIAL}\\

\vspace{3mm}

{\bf TAIL ESTIMATES FOR U-STATISTICS }\\

\vspace{4mm}

{\bf E. Ostrovsky, L.Sirota.}\\

\vspace{5mm}

Department of Mathematics, Bar-Ilan University, Ramat-Gan, 59200,Israel.\\
e-mail: \ eugostrovsky@list.ru \\

\vspace{3mm}

Department of Mathematics, Bar-Ilan University, Ramat-Gan, 59200,Israel.\\
e-mail: \ sirota3@bezeqint.net \\

\vspace{4mm}

{\bf Abstract.} \par

\vspace{3mm}

\end{center}

 \  We obtain in this paper a non-asymptotic  non-improvable up to multiplicative constant moment
and exponential tail estimates  for  distribution for $  U  \ - $ statistics by means of martingale representation. \par
 \ We show also the exactness of obtained estimations in one way or another by providing appropriate examples. \par

\vspace{4mm}

 \par

\vspace{3mm}


\vspace{4mm}

\section{ Introduction. Notations. Statement of problem.} \par

\vspace{3mm}

 \hspace{4mm}   Let $ (\Omega,F,{\bf P} ) $ be a probabilistic space, which will be presumed sufficiently
rich when we construct  examples (counterexamples).
 \ Let $ \{ \xi(i) \}, \ i = 1,2,\ldots, n, \ $ be independent identically distributed (i., i.d.)
random variables (r.v.) with values in the certain measurable space $  (X,S) , \  \Phi = \Phi(x(1), x(2), \ldots, x(d)) $
be a symmetric measurable non-trivial numerical function (kernel) of $  d  $ variables:
$  \Phi: X^d \to R, \ U(n) = U_n = U(n, \Phi,d) =  $

$$
 U(n, \Phi,d; \{  \xi(i) \} ) = {n \choose d}^{-1} \sum_{I \in I(d,n)}
\Phi( \xi(i_1), \xi(i_2), \ldots, \xi(i(d))), \ n > d \eqno(1.0)
$$
be a so-called $ U \ - $ statistic. Denote $ \deg  \Phi = d,  $

$$
\Phi = \Phi( \xi(1), \xi(2), \ldots, \xi(d)), \ r = \rank  \Phi \in [ 1,2, \ldots, d-1 ],\eqno(1.1)
$$

$$
\sigma(n) = \sigma_n = \sqrt{\Var (U(n))}, \hspace{4mm} \overline{T}(\Phi,x) := \sup_{ n > d} T(U(n) - {\bf E}U(n))/\sigma(n),x), \eqno(1.2)
$$
i.e. the uniform tail function for our $  U \ - $ statistics under natural norming. \par

 \ Let also $ I = I(n) = I(d; n) = \{i_1; i_2; \ldots; i_d \} $ be the set of indices of the
form $ I(n) = I(d; n) = \{ \vec{i} \}  = \{i \} = \{i_1, i_2, \ldots, i_d \} $ such that $ 1 \le i_1 < i_2 <
i_3 < i_{d-1} < i_d \le n; \  J = J(n) = J(d; n) $ be the set of indices of
the form (subset of  $ I(d; n))  \  J(d; n) = J(n) =   \{ \vec{ j} \} = \{ j  \} = \{j_1; j_2; \ldots; j_{d-1} \} $
such that $ 1 \le j_1 < j_2 \ldots < j_{d-1} \le n- 1. $ \par

\vspace{4mm}

 \ Recall that

$$
\sigma^2(n) = \Var(U(n)) \asymp n^{-r}, \ n \to \infty.
$$

\vspace{4mm}

  \ The {\it martingale representation } for the $ U \ - $ statistics as well as the
exact value for its variance $ \sigma^2(n) = \Var(U(n)) $ may be found, e.g. in \cite{Hoefding1}, \cite{Korolyuk1}, chapter 1.
Namely,

$$
U(n) - {\bf E} U(n) = \sum_{m=r}^d {d \choose m} U_{n, m},
$$
where

$$
 U_{n, m} = {n \choose m}^{-1} \sum \sum \ldots \sum_{1 \le i(1) < i(2) \ldots < i(m)  \le n}
 g_i( \xi(i(1)), \xi(i(2)),  \ldots, \xi(i(m))),
$$

$$
g_i( x(1), x(2),  \ldots, x(m)) =
$$

$$
\int_X \int_X \ldots \int_X \Phi(y(1), y(2), \ldots, y(m)) \prod_{s=1}^m
( \delta_{x(s)}(d y(s) ) - {\bf P}(d y(s) )) \times \prod_{s=m+1}^d {\bf P}(d y(s)).
$$

 \ The sequence

$$
 n \to S_m(n) = {n \choose m} U_{n,m}
$$
relative the natural filtration

$$
F_k = \sigma \{ \xi(1), \xi(2), \ldots, \xi(k) \},   \ F_0 = \{ \emptyset,X \}
$$
forme a martingale. \par

 \ Herewith

$$
\sigma^2(n) = \sum_{m=r}^d { d \choose m }^2 \cdot { n \choose m}^{-1} \cdot \Var \Phi =
$$

$$
r! \ {d \choose r}^2 \ n^{-r} \ \Var \Phi + O(n^{- r- 1}), \ n \to \infty.
$$

 \ Note in addition that it follows from Iensen inequality $  | \ g_i( \xi(i(1)), \xi(i(2)),  \ldots, \xi(i(m))) \ |_p \le | \ \Phi \ |_p. $

\vspace{4mm}

 \ Here and in the future for any r.v. $ \ \eta \ $  the function $ \ T_{\eta} ( x) \ $ will be denote its {\it tail} function:

$$
 T_{\eta}(x) \stackrel{def}{=} \max( {\bf P}(\eta > x), {\bf P}(\eta < - x) ), \ x > 0. \eqno(1.3)
$$

\vspace{3mm}

 \ We denote as usually the $ L(p) $ norm of the r.v. $ \eta $ as follows:

 $$
 |\eta|_p = \left[ {\bf E} |\eta|^p  \right]^{1/p}, \ p \ge 1;
 $$
and correspondingly

$$
M(p) = M(p, \Phi) =  M(p, \Phi, \{ \xi(i)  \}) \stackrel{def}{=} \sup_n | \ (U(n)  - {\bf E}U(n))/\sigma(n) \ |_p. \eqno(1.4)
$$

\vspace{4mm}

 {\bf We will derive the non-refined up to multiplicative constant
 moment and exponential tail estimations for distribution of normed  $  U \ - $ statistics,  \\
 indeed, to estimate the variables $ M(p, \Phi, \{ \xi(i)  \})  $ and} $ \overline{T}_{U(n)/\sigma(n)}(x). $ \par

\vspace{3mm}

  \ Evidently, these estimates may be applied for building  of a non-asymptotical confidence interval for unknown parameter
 by using the $ U \ -  $ statistics in the statistical estimation. \par

 \vspace{3mm}

 \ There are many works about this problem; the  next list is far from being complete:
  \cite{Borisov1}, \cite{Pena1}, \cite{Gine1}, \cite{Hoefding1}, \cite{Korolyuk1}, \cite{Kwapien1}, \cite{Ostrovsky10} etc.; see
 also reference therein. \par
 \ Notice that in the  classical book  \cite{Korolyuk1} there are many examples of applying of the theory of $ U \ -  $ statistics.
A new application, namely, in the modern adaptive estimation in the non - parametrical statistics may be found in the article
\cite{Bobrov1} and in the book \cite{Ostrovsky1}, chapter 5, section 5.13. \par

  \vspace{4mm}

\section{ Main result:  moments estimation for U \ - statistics. }

\vspace{4mm}

  \hspace{3mm} It is reasonable to suppose that $ {\bf E} \Phi = 0, \ \Var(\Phi) \in (0,\infty); $  moreover, we  can and will
assume without loss of generality $  \Var(\Phi) = 1, $   as long as it is constant;  and that all the moments of r.v. $ \Phi $
which are written below there exist; otherwise  it nothing to prove. \par

\vspace{4mm}

{\bf Theorem 2.1.} {\it Let  } $  {\bf E} \Phi = 0, \ \Var \Phi = 1, \ | \ \Phi \ |_p < \infty  $ {\it for some  value } $ p \ge 2. $
{\it Then  }

$$
\left| \frac{U(n)}{\sigma(n)} \right|_p \le C(d,r) \cdot  \left[\frac{p}{\log p} \right]^d \cdot |\Phi|_p, \ p \ge 2.\eqno(2.1)
$$

\vspace{4mm}

 {\bf Proof.} \par

 \vspace{3mm}

\ {\bf 1. Previous result.} The most recent (foregoing) results in this direction was obtained in \cite{Ostrovsky4}:

$$
\left| \frac{U(n)}{\sigma(n)} \right|_p \le C_0(d,r) \cdot  \left[\frac{p^d}{\log p} \right] \cdot |\Phi|_p, \ p \ge 2. \eqno(2.2)
$$

 \ Thus, we update a dependency on the  degree $  p. $    See also \cite{Hitczenko1}. \par

\vspace{3mm}

{\bf 2. Outline of the proof.} \ The inequality (2.2)  was obtained in \cite{Ostrovsky4} by means of the so-called
{\it martingale representation } for the
$  U \ - $ statistics, see \cite{Hoefding1}, \cite{Korolyuk1}, chapters 1,2; and using further the moment estimation for the centered
homogeneous of the degree $  d  $ polynomial martingales $ (\zeta_n, G_n), $ see \cite{Ostrovsky4}, of the form

$$
\sup_n \left|  \frac{\zeta_n}{\sqrt{\Var(\zeta_n)}} \right|_p \le C_2(d) \frac{p^d}{\log p}, \eqno(2.3)
$$
if of course $  | \ \zeta_n \ |_p < \infty. $ \par
 \ For the multiply series in rearrangement invariant spaces analogous result was obtained by S.V.Astashkin in
 \cite{Astashkin1}.  More information about martingale inequalities may be found  in the many works of
 D.L.Burkholder, see e.g. the articles  \cite{Burkholder1}-\cite{Burkholder3}. A famous
 survey on the martingale inequalities belongs to G.Peshkir and A.N.Shirjaev \cite{Peshkir1}. \par

 \ But in the recent publication about martingales \cite{Ostrovsky10} relaying in turn on the famous result  belonging to A.Osekowski
 \cite{Osekowski1} the estimate(2.4) was improved:

$$
\sup_n \left|  \frac{\zeta_n}{\sqrt{\Var(\zeta_n)}} \right|_p \le C_3(d) \left[\frac{p}{\log p}\right]^d. \eqno(2.4)
$$

 \ More exactly, the following important function was introduced by A.Osekowski (up to factor 2)
 in the article \cite{Osekowski1}:

$$
Os(p) \stackrel{def}{=} 4 \  \sqrt{2} \cdot \left( \frac{p}{4} + 1 \right)^{1/p} \cdot
\left( 1 + \frac{p}{\ln (p/2)}  \right), \ p \ge 4; \eqno(2.5)
$$
the case $ \ p \in [2, 4) \ $ is simple and may be considered separately.\par
 \ Note that

$$
K = K_{Os} \stackrel{def}{=} \sup_{p \ge 4} \left[\frac{Os(p)}{p/\ln p} \right] \approx 15.7858, \eqno(2.6)
$$
 is the so-called Osekowski's constant. \par
 \ Let us define the following numerical sequence $ \gamma(d), \ d = 1,2,\ldots:  \ \gamma(1) := K_{Os} = K,   $ (initial condition) and
by the following recursion

$$
\gamma(d+1) = \gamma(d) \cdot K_{Os} \cdot \left( 1 + \frac{1}{d} \right)^d. \eqno(2.7)
$$
 Since

$$
\left( 1 + \frac{1}{d} \right)^d \le e,
$$
we conclude

$$
\gamma(d) \le K_{Os}^d \cdot e^{d-1}, \ d = 1,2,\ldots. \eqno(2.8)
$$

 \ It is proved in \cite{Ostrovsky10} in particular that the "constant" $ C_3(d) $ in (2.4) allows the following simple estimate:
$   C_3(d) \le \gamma(d). $ \par

\vspace{3mm}

 \ The inequality (2.8) represents nothing more than $  d \ - $ dimensional and martingale generalization
of a classical Rosenthal's inequality for sums of independent random variables,  \cite{Rosenthal1}, see also \cite{Kallenberg1},
the exact values of constants in the Rosenthal's inequality see in \cite{Ostrovsky9}. \par

\ We apply further the more modern estimate (2.5) instead (2.4) into the considerations of the report \cite{Ostrovsky4}, we obtain
what is desired. \par

\vspace{3mm}

{\bf 3. Some details.}  {\it Let as before $  {\bf E} \Phi = 0, \ \Var \Phi = 1, \ | \ \Phi \ |_p < \infty $ for some  value $ p \ge 2. $
 Let also the sequence   } $ \gamma(d) $ {\it  be defined in (2.7) (and in (2.8)). Then  }

 \vspace{3mm}

$$
\left|U(n) \right|_p \le \sum_{m=r}^d \gamma(m) \cdot { d \choose m } \cdot { n \choose m }^{-1/2} \cdot
 \left( \frac{p}{\ln p} \right)^m \cdot | \ \Phi \  |_p, \ p \ge 2. \eqno(2.9)
$$

{\it and in turn after evident simplification}

$$
| \ U(n) \ |_p  \le C(d,r) \ n^{-r/2} \ \left[\frac{p}{\ln p} \right]^d \ | \ \Phi \ |_p, \ p \ge 2. \eqno(2.10)
$$

\vspace{3mm}

 {\bf  Proof}.  \ We can write using the martingale representation for $ \ U \ - \ $ statistics

$$
U(n) = \sum_{m=r}^d { d \choose m} \ { n \choose m }^{-1} \ \sum_{j \in J(m,n)} \mu_j, \eqno(2.11)
$$
where $ (\mu_k, F_k) $ is certain centered martingale,

$$
 \card J(m,n) = { n \choose m }, \  \Var \mu_j = \Var \Phi = 1.
$$
 Denote for brevity

 $$
 N = N(m,n) = { n \choose m },
 $$
  then

$$
U(n) = \sum_{m=r}^d { d \choose m}  \ N^{-1/2} (m,n) \ \zeta(m,n). \eqno(2.12)
$$

 \ We apply the triangle inequality for the $  L(p) $ norm:

$$
| \ U(n) \ |_p \le  \sum_{m=r}^d { d \choose m}  \ N^{-1/2} (m,n) \ | \ \zeta(m,n) \ |_p. \eqno(2.13)
$$

\ Each term $ \zeta(m,n) $ is the centered  polynomial martingale of degree $  m   $ generated by the function $ \Phi. $
One can apply the Osekowski's inequality (2.4):

$$
| \ \zeta(m,n) \ |_p \le \gamma(m) \ \left[ \frac{p}{\ln p} \right]^m \ | \ \Phi \ |_p. \eqno(2.14)
$$
 \ It remains to substitute into (2.13). \par

\vspace{3mm}

{\bf 4.}  The assertion of theorem 2.1 follows immediately from the estimate (2.9), but this estimate gives us certain numerical
estimate for the value $ \ |U_n|_p.  $ \par

\vspace{3mm}

{\bf Remark 2.1.} As long as $  | \ U_n \ |_p \ge   | \ U_n \ |_2 = \sigma_n, $ we deduce that
every time when  $  | \ \Phi \ |_p < \infty,  $

$$
| \ U(n) \ |_p \asymp n^{-r/2} \asymp \sigma(n), \ n \to \infty.  \eqno(2.15)
$$

 \ We generalize further  this relation on more general than $ L_p(\Omega) $ spaces.\par

\vspace{4mm}

 \ Let us discuss now the {\it lower bounds} for the theorem 2.1. To be more precise, we denote

$$
K_d(p) \stackrel{def}{=} \sup_{ \{ \xi(i)  \}  }
\sup_{ 0 \ne \Phi \in L_p } \sup_n  \left\{ \frac{ M(p, \Phi, \{ \xi(i)  \})}{ n^{-r/2} \ |\Phi|_p} \right\}. \eqno(2.16)
$$

\vspace{4mm}

{\bf Theorem 2.2}. {\it We conclude taking into account formulated above our definition and restrictions}

$$
K_d(p) \asymp \left[ \frac{p}{\ln p} \right]^d, \ p \ge 2. \eqno(2.17)
$$

\vspace{4mm}

 \ {\bf Proof.} It remains to prove only the lower bound for the value $ K_d(p). $ \par
 \ The relation (2.12)  has been conjectured (hypothesis) in the article  \cite{Ostrovsky4}, page 19 for the polynomial martingales.
It was proved (before!) for the polynomial martingales from appropriate independent random variables in \cite{Hitczenko1}, \cite{Kallenberg1}.
The case of arbitrary polynomial martingales was grounded in authors report \cite{Ostrovsky10}. \par
 \ We represent here a very simple example in order to obtain the bottom border for $ K_d(p) $ exactly for the $  U \ - $ statistics
 still for arbitrary value $ d = 1,2,\ldots. $ Let $  n = 1 $  and let a r.v. $ \eta $ has a standard Poisson distribution with unit parameter

$$
{\bf P}(\eta = k) = e^{-1}/k!, \ k = 0,1,2,\ldots
$$
and define $ \xi = \eta - 1; $ then $ \xi  $  is centered, $  \Var (\xi) = 1  $ and it is no hard to calculate

$$
|\xi|_p \sim \frac{p}{e \cdot \ln p}, \ p \to \infty.
$$

 \ Let also $  \xi(i) $ be independent copies of $ \xi. $   Then

$$
| \ \prod_{i=1}^d \xi(i) \ |_p \sim e^{-d} \left[ \frac{p}{\ln p} \right]^d, \ p \to \infty. \eqno(2.18)
$$
 Therefore

$$
\underline{\lim}_{p \to \infty} \left\{ \ K_d(p) : \left[ \frac{p}{\ln p} \ \right]^d \  \right\} \ge e^{-d}. \eqno(2.19)
$$

 \vspace{3mm}

 \ Thus, obtained in this section estimate (2.1) of theorem (2.1) is essentially non-improvable relative the parameter $  p $
for all the values of dimension $  d, $ of course, up to multiplicative constant. \par

 \vspace{4mm}

\section{ Estimations of U-statistics in the Grand Lebesgue Spaces and in the exponential Orlicz space norms. }

\vspace{4mm}

 \ Let $  \psi = \psi(p), \ p \in [2,b), \ b = \const, \ 1 < b \le \infty   $  (or $ p \in [1,b] ) $  be certain bounded from below:
 $ \inf \psi(p) > 1 $  continuous
inside the {\it  semi-open  } interval $ [2, b) $ numerical function. We can and will suppose

$$
b = \sup \{p, \ \psi(p) < \infty \},
$$
so that $ \supp \psi = [2,b)  $ or $ \supp \psi = [2,b].$ The set of all such a functions will be denoted by
$ \Psi(b); \ \Psi:= \Psi(\infty). $ \par
 \ For each such a function $  \psi \in \Psi(b) $ we define

$$
\psi_d(p) \stackrel{def}{=}  \left[ \frac{p}{ \ln p} \right]^d \cdot \psi(p). \eqno(3.1)
$$
 \ Evidently,  $ \psi_d(\cdot) \in \Psi(b). $ \par

 \ By definition, the (Banach) space $  G\psi = G\psi(b)  $  consists on all the numerical valued random variables $ \{ \zeta \}  $
defined on our  probability space $ (\Omega,F,{\bf P} ) $ and having a finite norm

$$
||\zeta|| G\psi \stackrel{def}{=} \sup_{p \in (1,b)} \left[ \ \frac{|\zeta|_p}{\psi(p)} \ \right] < \infty. \eqno(3.2)
$$
 \ These spaces are suitable in particular for an investigation  of the random variables and the random processes (fields) with
exponential decreasing tails of distributions, the Central Limit Theorem in Banach spaces, study of Partial Differential
Equations etc., see e.g. \cite{Kahane1}, \cite{Buldygin1}, \cite{Kozatchenko1}, \cite{Ostrovsky1}, chapter 1,
\cite{Fiorenza1}-\cite{Fiorenza3}, \cite{Iwaniec1}-\cite{Iwaniec2} etc. \par
 \ More detail, suppose $  0 < ||\zeta|| := ||\zeta||G\psi < \infty. $ Define the function

$$
\nu(p) = \nu_{\psi}(p)  = p \ln \psi(p), \ 2 \le p < b
$$
and put formally $ \nu(p) := \infty, \ p < 2  $ or $ p > b. $ Recall that the Young-Fenchel, or Legendre transform $ f^*(y) $
for arbitrary function $  f: R \to R $ is defined (in the one-dimensional case) as follows

$$
f^*(y) \stackrel{def}{=} \sup_x (x y - f(x)).
$$

 \ It is known that

 $$
 T_{\zeta}(y) \le \exp \left( - \nu_{\psi}^*(\ln (y/||\zeta||) )  \right), \ y > e \cdot ||\zeta||. \eqno(3.3)
 $$
 \ Conversely, if (3.3) there holds in the following version:

$$
T_{\zeta}(y) \le \exp \left( - \nu_{\psi}^*(\ln (y/K) )  \right), \ y > e \cdot K, \ K = \const > 0, \eqno(3.4)
$$
and the function $ \nu_{\zeta}(p), \ 2 \le p < \infty $ is positive, continuous, convex and such that

$$
\lim_{p \to \infty} \psi(p) = \infty,
$$
then $ \zeta \in G\psi $ and besides

$$
||\zeta||G\psi \le C(\psi) \cdot K. \eqno(3.5)
$$

\ Moreover, let us introduce the {\it exponential }  Orlicz space $ L^{(M)} $ over the source probability space
$ (\Omega,F,{\bf P} ) $ with proper Young-Orlicz function

$$
M(u) := \exp \left(  \nu_{\psi}^*(\ln |u| )  \right), \ |u| > e
$$
or correspondingly

$$
M_d(u) := \exp \left(  \nu_{\psi_d}^*(\ln |u| )  \right), \ |u| > e
$$
and as ordinary $  M(u) = M_d(u) = \exp(C \ u^2) - 1, \ |u| \le e. $ It is known \ \cite{Ostrovsky11} that the $ G\psi $ norm
of arbitrary r.v. $  \zeta  $ is complete equivalent to the its norm in Orlicz space $ L^{(M)}: $

$$
||\zeta||G\psi \le C_1 ||\zeta||L^{(M)} \le C_2||\zeta||G\psi, \ 1 \le C_1 \le C_2 < \infty;
$$

$$
||\zeta||G\psi_d \le C_3 ||\zeta||L^{(M_d)} \le C_4||\zeta||G\psi_d, \ 1 \le C_3 \le C_4 < \infty.
$$

\vspace{4mm}

{\bf Example 3.1.} The estimate for the r.v. $ \xi $ of a form

$$
|\xi|_p  \le C_1 \ p^{1/m} \ \ln^r p, \ p \ge 2,
$$
where $ C_1 = \const > 0, \ m = \const > 0, \ r = \const, \ $ is quite equivalent to the following tail estimate

$$
T_{\xi}(x) \le  \exp \left\{ - C_2(C_1,m,r) \ x^m \ \log^{-m r}x \right\}, \ x > e.
$$

\vspace{3mm}

 \ It is important to note that the inequality (3.4) may be applied still when the r.v. $ \xi $ does not have the
exponential  moment, i.e. does not satisfy the famous Kramer's condition. Namely, let us consider next example. \\

\vspace{3mm}

 {\bf Example 3.2.} Define the following $ \Psi \ -  $ function.

$$
\psi_{[\beta]}(p) := \exp \left( C_3 \ p^{\beta} \right), \  p \in [2, \infty), \ \beta = \const > 0.
$$

 \ The r.v. $  \xi $ belongs to the space $ G \psi_{[\beta]} $ if and only if

$$
T_{\xi}(x) \le \exp \left( - C_4(C_3,\beta) \ [\ln (1 + x)]^{1 + 1/\beta}   \right), \ x \ge 0.
$$

 \ See also \cite{Ostrovsky4}. \par

\vspace{3mm}

 \ Let us return to the source problem. Assume that there exists certain function $  \psi(\cdot) \in \Psi(b), \   b = \const \in (2, \infty)  $
such that $ \Phi \in G\psi(b). $ For instance, this function may be picked by the following {\it natural} way:

$$
\psi_{\Phi}(p) := |\Phi|_p, \eqno(3.6)
$$
if of course there exists and is finite at last for some value $  p  $ greatest than 2,  obviously,
with the correspondent value $  b. $\par

\vspace{3mm}

{\bf Theorem 3.1.} We propose under formulated above conditions

$$
\sup_n || U(n)/\sigma(n)||G\psi_d \le  C(\psi) \ ||\Phi||G\psi. \eqno(3.7)
$$

\vspace{3mm}

{\bf Proof} is very simple. Assume $ ||\Phi||G\psi  \in (0, \infty). $ It follows immediately from the direct definition
of the norm in the Grand Lebesgue Spaces

$$
| \ \Phi \ |_p \le ||\Phi||G\psi \cdot \psi(p).
$$
 \ We apply the inequality (2.1) of theorem 2.1 for the values $  p  $ from the set $  p \in [2,b): $

$$
\sup_n | \ U(n)/\sigma(n) \ |_p \le C(d,r) \cdot  \left[\frac{p}{\log p} \right]^d \cdot |\Phi|_p \le
$$

$$
C(d,r) \cdot  \left[\frac{p}{\log p} \right]^d \cdot ||\Phi||G\psi \cdot \psi(p) =
$$

$$
C(d,r) \cdot  ||\Phi||G\psi \cdot \psi_d(p),
$$
or equally

$$
\sup_n || \ U(n)/\sigma(n) \ ||G\psi_d \le C(d,r) \cdot  ||\Phi||G\psi.
$$

\vspace{3mm}

{\bf Example 3.3.} Suppose

$$
T_{\Phi}(x) \le  \exp \left\{ - C_5 \ x^m \ \log^{-m r}x \right\}, \ x > e, \eqno(3.8)
$$
where $ C_5 = \const > 0, \ m = \const > 0, \ r = \const. \ $   As we know,

$$
|\Phi|_p  \le C_6(C_5,m,r) \ p^{1/m} \ \ln^r p, \ p \ge 2.
$$

 \ We use theorem 3.1

$$
\sup_n | \ U(n)/\sigma(n) \ |_p \le C_7 \ p^{d + 1/m} \ (\ln p)^{ r - d  }, p \ge 2,
$$
and  we conclude returning to the tail of distribution

$$
\sup_n T_{U(n)/\sigma(n)}(x) \le  \exp \left\{ - C_8 \ x^{m/(1 + dm)} \ \log^{-m (r-d)/(1 + dm)}x \ \right\}, \ x \ge e. \eqno(3.9)
$$

 \ Thus, we obtained in this way the {\it exponential } bounds for distribution for the normed $  U \ -  $ statistics,
expressed only through the very simple source data (3.8). \par

 But the authors are not convinced of the finality of these estimations in the considered case; cf., e.g.
\cite{Borisov1}, \cite{Gine1}, \cite{Korolyuk1}, chapter 2, \cite{Pena1}.\par

\vspace{3mm}

{\bf Example 3.4.} Assume now

$$
T_{\Phi}(x) \le \exp \left( - C_9 \ [\ln (1 + x)]^{1 + 1/\beta}   \right), \ x \ge 0.\eqno(3.10)
$$
or more strictly

$$
\exp \left( - \tilde{C}_9 \ [\ln (1 + x)]^{1 + 1/\beta}   \right) \le   T_{\Phi}(x) \le
$$

$$
\exp \left( - C_9 \ [\ln (1 + x)]^{1 + 1/\beta}   \right), \ x \ge 0,   0 < C_9 \le \tilde{C}_9 < \infty.  \eqno(3.10a)
$$

 \ Then

$$
|\ \Phi \ |_p \le \exp \left( C_{10} \ p^{\beta} \right),   \ p \ge 2,
$$
therefore by virtue of theorems 2.1 and  3.1

$$
\sup_n | \ U(n)/\sigma(n) \ |_p  \le C_{11} \ p^d \  (\ln p)^{-d} \ \exp \left( C_{10} \ p^{\beta} \right) \le
$$

$$
\exp \left( C_{12} \ p^{\beta} \right), \ p \ge 2,
$$
and we have returning again to the tails of distribution

$$
\sup_n T_{U(n)/\sigma(n)}(x) \le \exp \left( - C_{13} \ [\ln (1 + x)]^{1 + 1/\beta}   \right), \ x \ge 0.\eqno(3.11)
$$

 \ As long as

$$
\sup_n T_{U(n)/\sigma(n)}(x) \ge T_{U(1)/\sigma(1)}(x) =  T_{\Phi}(x),
$$
cf. (3.10) and (3.10a), we conclude that the assertion of theorem 3.1, i.e. exponential bound for tail distribution,
is also in general case non improvable. \par

\vspace{4mm}

{\bf Remark 3.1.} As we promised to prove, the proposition of Remark 2.1 it remains valid still for the Grand Lebesgue Spaces
$  G\psi $  and for exponential Orlicz spaces of the form $ L^{(M)}: $
 every time when  $  || \ \Phi \ ||G\psi < \infty  $ or equally $  || \ \Phi \ ||L^{(M)} < \infty  $

$$
|| \ U(n) \ ||G\psi_d \asymp n^{-r/2} \asymp \sigma(n) \asymp || \ U(n) \ ||L^{(M_d)},  \ n \to \infty. \eqno(3.12)
$$

 \vspace{4mm}

\section{Concluding remarks. }

\vspace{4mm}

 \hspace{4mm}  {\bf A.} It is interest, by our opinion, to obtain  analogous estimates for dependent source random variables,
for instance, for martingales or mixingales. Some preliminary results in this directions may be found in \cite{Borisov1},
 \cite{Gine1}. \\

 \vspace{3mm}

{\bf B.} We do not aim to derive the best possible  values of appeared in this report constants.
It remains to be done. \\

\vspace{3mm}

{\bf C.} Perhaps, the case of the so-called  $ V \ - $ statistics  may be investigated analogously. \par


\vspace{4mm}

\begin{thebibliography}{99}


\bibitem{Astashkin1}
{\sc Astashkin S.V.} (1999). {\it Multiply Series in Rearrangement Invariant Spaces. }
Funct. Anal. Applic.,  $ N^o $  2, V. 33, 141 - 143.

\bibitem{Bobrov1}
{\sc Bobrov P.B., Ostrovsky E.I.} (1997). {\it The Adaptive Estimation of
a Regression, Density and Spectrum.} Probability and Statistics, Collective Works
of St. Petersburg Branch of the Steklov Math. Institute, (in Russian), (2), v. 244,
28 \ - \ 45.

\bibitem{Borisov1}
{\sc I. S. Borisov, N. V. Volodko.} {\it A note on exponential inequalities
for the distribution tails of canonical Von Mises statistics of dependent observations.}\\
arXiv:0907.0058v3 [math.PR] 27 Feb 2015

\bibitem{Buldygin1}
{\sc Buldygin V.V., Kozachenko Yu.V.} (2000). {\it Metric Characterization
of Random Variables and Random Processes. } AMS, 678, Providence, R.I.

\bibitem{Burkholder1}
{\sc D.L.Burkholder.} {\it Distribution Functions Inequalities for Martingales.} Ann. Probab., {\bf 1}, (1973),
19-42.

\bibitem{Burkholder2}
{\sc D.L.Burkholder.} {\it Explorations in Martingale Theory and its Applications.} Ecole d'Ete de
Probabilities de Saint-Flour XIX, 1989, pp. 1-66, Lecture Notes in Math., 1464, Springer, Berlin, 1991.

\bibitem{Burkholder3}
{\sc  D.L.Burkholder.} (1998). {\it Sharp Inequalities for Martingales and
Stochastic Integrals.} Asterisque, 157/158; 75-96.

\bibitem{Choi1}
{\sc Choi K.P.} {\it Some sharp inequalities for martingale transform.} Trans. Amer.
Math. Soc., 307, No 1, (1988), 279-300, MR0936817.

\bibitem{Fiorenza1}
   {\sc Capone C., Fiorenza A., Krbec M.} {\it On the Extrapolation Blowups in the
   $ L_p $ Scale. } Collectanea Mathematica, {\bf 48}, 2, (1998), 71-88.

 \bibitem{Fiorenza2}
 {\sc Fiorenza A.} {\it Duality and reflexivity in grand Lebesgue spaces.}
       Collectanea Mathematica (electronic version), {\bf 51}, 2, (2000), 131-148.

\bibitem{Fiorenza3}
 {\sc Fiorenza A., and Karadzhov G.E.} {\it Grand and small Lebesgue spaces and
       their analogs. } Consiglio Nationale Delle Ricerche, Instituto per le
      Applicazioni del Calcoto Mauro Picone", Sezione di Napoli, Rapporto tecnico n.
      272/03, (2005).

\bibitem{Gine1}
{\sc Gine, E., R. Latala, and J. Zinn.} (2000). {\it Exponential and moment inequalities for U-statistics.}

\bibitem{Hitczenko1}
{\sc P. Hitczenko.} {\it Best constants in martingale version of Rosenthal inequality.} Ann. Probab. 18
No. 4 (1990), 1656-1668.

\bibitem{Hoefding1}
{\sc Hoefding, W.} (1948). {\it A class of statistics with asymptotically normal distribution.} The
Annals of Mathematical Statistics, V.1,  293-325.

\bibitem{Iwaniec1}
   {\sc Iwaniec T., and Sbordone C.} {\it On the integrability of the Jacobian under
      minimal hypotheses. } Arch. Rat.Mech. Anal., 119, (1992), 129–143.

\bibitem{Iwaniec2}
 {\sc Iwaniec T., P. Koskela P., and Onninen J.} {\it Mapping of finite distortion:
   Monotonicity and Continuity.}  Invent. Math. 144 (2001), 507-531.

\bibitem{Kahane1}
{\sc Kahane J.P.} {\it Properties locales des fonctions a series de Fourier aleatoires.}
 Studia Math. (1960), 19, No 1, 1-25.

\bibitem{Kallenberg1}
{\sc Kallenberg Olav and Sztencel Rafal.} {\it Some dimension-free features
of vector-valued martingales.} Probability Theory Related Fields.
1991, {\bf 88,}  215-247.

\bibitem{Korolyuk1}
{\sc Korolyuk V.S., Borovskokh Yu.V. }  (1994). {\it Theory of U-Statistics.}
Kluwer Verlag, Dodrecht.

\bibitem{Kozatchenko1}
{\sc Kozatchenko Yu. V., Ostrovsky E.I.} {\it The Banach Spaces of
random Variables of subgaussian type.} Theory of Probab. and Math.
Stat. (in Russian). Kiev, KSU, 32, 43-57, (1985).

\bibitem{Kwapien1}
{\sc  Kwapien Stanislaw.} {\it  On Hoefding decomposition in L(p).}
Illinois Journal of Mathematics, Volume 54, Number 3, 2010, Pages 1205-1211.

\bibitem{Osekowski1}
{\sc Osekowski A.} {\it A Note on Burkholder-Rosenthal Inequality.  }
Bull. Polish Academy of Science, Math., {\bf 60}, (2012), 177-185.

\bibitem{Ostrovsky1}
{\sc Ostrovsky E.I.} {\it Exponential estimations for Random Fields and
     its applications.} (in Russian). 1999, Moskow-Obninsk, Russia, OINPE.

\bibitem{Ostrovsky4}
{\sc  Ostrovsky E. } {\it Bide-side exponential and moment inequalities  for
      tails of distribution of Polynomial Martingales.}\\
      arXiv: math.PR/0406532 v.1  Jun. 2004

\bibitem{Ostrovsky7}
{\sc Ostrovsky E. and  Sirota L. } {\it Moment and tail inequalities for polynomial martingales.
 The case of heavy tails.} \\
arXiv:1112.2768v1 [math.PR] 13 Dez 2011

\bibitem{Ostrovsky9}
{\sc Ostrovsky E. and  Sirota L.} {\it Schl\"omilch and Bell series for Bessel's functions, with probabilistic applications.}\\
arXiv:0804.0089v1 [math.CV] 1 Apr 2008

\bibitem{Ostrovsky10}
{\sc Ostrovsky E. and  Sirota L.} {\it  Sharp moment estimations for polynomial martingales.} \\
arXiv:1410.0739v1 [math.PR] 3 Oct 2014

\bibitem{Ostrovsky11}
{\sc Ostrovsky E. and  Sirota L.} {\it  Banach spaces characterization of random vectiors with exponential
decreasing tails of distribution. } \\
arXiv:1601.04766v1 [math.PR] 19 Jan 2016

\bibitem{Pena1}
{\sc De la Pen'a, V. H.} (1992). {\it Decoupling and Khintchine's inequalities for U-statistics.} The
Annals of Probability, V.4, 1877-1892.

\bibitem{Peshkir1}
{\sc Peshkir G., Shirjaev A.N.} {\it The Khintchine inequalities and martingale
expanding sphere of their action.} Russian Math. Surveys; {\bf 50,} 5, 849-904, (1995).

\bibitem{Rosenthal1}
{\sc Rosenthal H.P.}  {\it On the Subspaces of $ L_p, \ (p > 2) $ spanned by Sequences of
independent Variables. } Israel J. Math., 1970, V.3, pp. 253-273.

\end{thebibliography}
\end{document}